\documentclass{article}
\usepackage{amsmath,amssymb,url,latexsym,amsthm}
\usepackage[all]{xy}
\usepackage{color}

\newtheorem{defi}{Definition}[section]
\newtheorem{theo}{Theorem}[section]
\newtheorem{lemma}{Lemma}[section]
\newtheorem{prop}{Proposition}[section]
\newtheorem{coro}{Corollary}[section]

\def\into{ \rightarrowtail }

\def\splito{ \rightleftarrows }
\def\trio{ \triangleright}

\newdir{ >}{{}*!/-8pt/@{>}}
\def\EE{ \mathbb{E} }

\def\VV{ \mathbb{V} }

\def\into{ \rightarrowtail }

\def\splito{ \rightleftarrows }

\newcommand{\Gp}{\mathsf{Gp}}
\newcommand{\Ab}{\mathsf{Ab}}

\newcommand{\Set}{\mathsf{Set}}
\newcommand{\Mag}{\mathsf{Mag}}
\newcommand{\Mon}{\mathsf{Mon}}

\newcommand{\Mal}{\mathsf{Mal}}
\newcommand{\Aff}{\mathsf{Aff}}
\newcommand{\BooRg}{\mathsf{BooRg}}

\newcommand{\CoM}{\mathsf{CoM}}
\newcommand{\Vect}{\mathsf{Vect}}
\newcommand{\Imp}{\mathsf{Imp}}
\newcommand{\HsLat}{\mathsf{HsLat}}
\newcommand{\Heyt}{\mathsf{Heyt}}

\newcommand{\Subt}{\mathsf{Subt}}

\def\trio{ \triangleright}

\begin{document}

\author{Dominique Bourn}

\title{On the concept of \emph{Algebraic Crystallography}}

\date{}

\maketitle

\begin{abstract}
 Category Theory provides us with a clear notion of what is an internal structure. This will allow us to focus our attention on a certain type of relationship between \emph{context} and 
 \emph{structure}.
\end{abstract}

\section{Introduction}\label{sec-Introduction}

It is known that, in the unital and strongly unital categories \cite{1} and in the subtractive varieties \cite{Urs} or categories \cite{ZJ} as well, on any object $X$ there is at most one structure of abelian group object. But, in these contexts, this did not seem so surprising because the four cases were closely related (strongly unital=unital+subtractive \cite{ZJ}) and because of the kind of their varietal origins: using an Eckmann-Hilton argument \cite{EH}, this uniqueness property arised naturally because, and when, some term in the definition of the varietal examples in question became a homomorphism in this variety.

Similar situations for other algebraic structures were even known from a long time; for instance, it was clear that in a pointed J\'onnson-Tarski variety, on any algebra $X$ there is at most one internal commutative  monoid structure; the same property holds for the commutative  and associative (=autonomous) Mal'tsev operations in the  Mal'tsev varieties \cite{Smith}. And again the limpid varietal contexts supplied the same simple explanation for this phenomenon.

But recently we were led to observe that the uniqueness structure for abelian group objects  still holds in the new context of Congruence Hyperextensible categories \cite{2}. This, in restrospect, emphasized that the uniqueness of the autonomous Mal'tsev operations was already noticed in the Congruence Modular Varieties \cite{Gumm} as well.

This phenomenon of uniqueness of some kinds of algebraic structure being now clearly extended to much larger contexts than the ones of the first paragraph, and the explanation by the existence of some kinds of terms in the definition of the varieties being no longer valid, it cannot remain possible to accept this uniqueness so easily and to keep it as an unquestioned property. So, we propose the following:
\begin{defi}
	A finitely complete category $\mathbb E$ is called crystallographic with respect to a given algebraic structure 
	$S$ when, on any object $X$ in $\mathbb E$, there is at most one internal algebraic structure of this kind. 	
\end{defi} 
This terminology is chosen because, in such a category $\mathbb E$, the algebraic structure $S$ in question is growing so punctually scarce. There are two extremal cases:\\ i) when on any object $X$ in $\EE$ there is one and only one $S$-structure, we say that $\EE$ is \emph{intensively crystallographic} with respect to the structure $S$;\\
ii) when the only $S$-structures are the terminal object $1$ and its subobjects, we say that $\EE$ \emph{trivializes} the structure $S$.\\
Concerning the structure of group, any additive category is intensively chrystallographic with respect to it, and any Mal'tsev variety whose Mal'tsev term $p$ satisfies the Pixley axiom $p(x,y,x)=x$ \cite{Pix} trivializes it.

Finally, when a structure $\Sigma$ is such that any category $\Sigma(\EE)$ of internal $\Sigma$-structures in $\EE$ trivializes the structure $S$, we say that the structure $\Sigma$ trivializes the structure $S$.

The aim of this work is to produce examples and to establish the very first properties and general questionings about this notion of \emph{Algebraic Crystallography}. It will lead to an unexpected and rather spectacular outcome with an example of a variety $\mathbb{H}$ which is crystallographic for the structure of abelian group and in which the punctual sarcity of abelian group objects seems to be offset by a kind of multiplication of this structure inside it. Indeed, the category $\Ab\mathbb{H}$ of abelian objects in this variety $\mathbb{H}$ will appear:\\ 
i) to fully faithfully embed the category $\Ab$ of abelian groups by a functor $h: \Ab \rightarrow \Ab\mathbb{H}$\\ and in an independant way\\
ii) to faithfully contain any category $K$-$\Vect$ of $K$-vector spaces, provided that the field $K$ is not of characteristic $2$, by a functor $w_{K}: K$-$\Vect\rightarrow \Ab\mathbb{H}$.\\
So, the abelian group  $(V,+)$ underlying a $K$-vector space $V$ will be represented by two distinct  objects $h(V,+)$ and $w_K(V,+)$ in the abelian category  $\Ab\mathbb{H}$.\\
We give an example of the same kind of multiplication process in a non-pointed context.

\section{First observations from the unital setting}\label{unital}

In this article any category  will be suppose finitely complete. The kernel equivalence relation of a map $f: X\to Y$ is denoted by $R[f]$. Given any algebraic structure $S$ and any category $\EE$, we denote by $S(\EE)$ the category of internal $S$-structures in $\EE$ and by $U^S_{\EE}: S(\EE)\to \EE$ the canonical forgetful functor. The category $S(\EE)$ is finitely complete as soon as so is $\EE$, and its terminal object is the unique possible $S$-structure on the terminal object $1$ of $\EE$. The functor $U^S_{\EE}$ is clearly left exact and it reflects isomorphims. 

Suppose now that the $S$-structure in question has no constant. Since, given any subobject $J \into 1$ of $1$ in $\EE$, we get $J^n\simeq J$ for any $1 \leq n$ in a unique possible way, there is one and only one $S$-structure on this object $J$. Whence a canonical fully faithful embedding $J_{\EE}: Sub_1\EE \into S(\EE)$ from the fully faithful subcategory $Sub_1\EE$ of the subobjects of $1$ in $\EE$ such that $U^S_{\EE}.J_{\EE}: Sub_1\EE \to \EE$ is nothing but the inclusion $Sub_1\EE \into \EE$. 

\subsection{General unital setting}

A \emph{unital} category \cite{1} is a pointed category $\mathbb E$ such that the canonical pair of injections: $$X\stackrel{l_X}{\rightarrowtail} X\times Y \stackrel{r_Y}{\leftarrowtail} Y \;\;\;\;\;\;\;\;\;\;\;\;\;\;\;\;\;\;\;\;\;\;\;\; (1)$$ is jointly strongly epic, i.e. such that the only (up to isomorphism) subobject of $X\times Y$ containing $l_X$ and $r_Y$ is $1_{X\times Y}$. 

A pointed variety $\VV$ is unital if and only if it is a  J\'onnson-Tarski variety \cite{BB}, i.e. if it has a unique constant $0$ and a binary term $+$ satisfying $x+0=x=0+x$. So, the main examples of unital category are the categories $\Mag$ and $\Mon$ of unitary magmas and monoids. The unital setting appears to be the right one in which there is an \emph{intrinsic notion of commutative pairs of subobjects}: given  any pair $u: U\rightarrowtail X$, $v: V\rightarrowtail X$ of subobjects, there is at most one factorization $\varphi$:
$$
\xymatrix@=10pt{
	&  U \ar@{ >->}[dl]_{l_U}  \ar@{ >->}[dr]^{u}  & & \\
	U \times V \ar@{.>}[rr]^{\varphi} &&  X   \\
	& V \ar@{ >->}[ul]^{r_V}  \ar@{ >->}[ur]_{v}  & & 
}
$$
 making the previous diagram commute. So, when such a map does exist, \emph{the subobjects $u$ and $v$ are said to commute} and the map $\varphi$ is called the {\em cooperator} of the pair. This situation is denoted by $[u,v]=0$. We can then follow the usual variations:
a subobject $u:U\into X$ is central when $[u,1_{X}]=0$;
an object $X$ is commutative  when $[1_{X},1_{X}]=0$.

By its cooperator $\varphi: X\times X \to X$, any commutative object $X$ is endowed with a  structure of internal unitary magma \emph{which turns out to be an internal commutative monoid}. Moreover, any morphism $f :X\to Y$ between two commutative objects preserves these monoid structures. When this commutative monoid is an internal abelian group, the object $X$ is said to be \emph{abelian}. 	Whence, in this context, the two natural intrinsic fully faithful subcategories: $\Ab(\EE)\hookrightarrow\CoM(\EE)\hookrightarrow \EE$ of the commutative and  abelian objects.

A collateral aspect of these observations was noticed:	
\emph{in a unital setting,
on an object $X$, there is atmost one structure of internal unitary magma}, 
but it was  not  thorough. We can now assert:
\begin{prop}
A unital category $\EE$ is crystallographic with respect to the structure of unitary magma and, a fortiori, to the structure of commutative monoid and abelian group.	
\end{prop}

\subsection{Special unital settings}

There are two extremal cases of unital categories: \cite{1}:
i) the only commutative object is the terminal object $1$;
ii) any object $X$ is commutative.

The second case holds when the diagram $(1)$ is actually the diagram of the sum of $X$ and $Y$; the pointed category $\EE$ is then called \emph{linear}. The first case holds when we have: $[u,v]=0 \implies u\cap v=0_X: 1\into X$; the pointed category $\EE$ is then said \emph{stiffly unital}.
These two extremal situations give rise to two extremal cases of crystallography:
\begin{defi}
A category $\EE$ trivializes (resp. strongly trivializes) a structure $S$ when the subobjects $J$ of $1$ are the only ones to have a $S$-structure (resp. the terminal object $1$ is the only one to have a $S$-structure).

A category $\EE$ is intensively chrystallographic with respect to a  structure $S$, when, on any object $X$, there is one and only one internal structure $S$.
\end{defi} 
Accordingly, the category $\EE$ (resp. strongly) trivializes the structure $S$ if and only if the canonical fully faithful embedding $J_{\EE}: Sub_1\EE \into S(\EE)$ is an equivalence of categories, where  $Sub_1\EE$ is the fully faithful subcategory of the subobjects of  $1$ in $\EE$ (resp. if and only if the canonical embedding $1 \into S(\EE)$ is an equivalence of categories where $1$ is the discrete category with only one object). And the  category $\EE$ is intensively chrystallographic with respect to the structure $S$ if and only if the forgetful functor $U^S_{\EE}: S(\EE) \to \EE$ is an equivalence of categories.

Any linear category is intensively chrystallogaphic with respect to the structure of monoid. The category $\BooRg$ of boolean rings and the category $\HsLat$ of Heyting semi-lattices  strongly trivialize the structure of  monoid, see \cite{BB}; this is the case as well for the dual $(\Set_*)^{op}$ of the category of pointed sets and more generally for the dual $(\EE_*)^{op}$ of the category of pointed objects in $\EE$ when it is topos, see \cite{1}. 
Finally we come to the following:
\begin{defi}
An algebraic structure $\Sigma$ (resp. strongly) trivializes another algebraic structure $S$  when, for any category $\EE$, the category $\Sigma(\EE)$ of internal  $\Sigma$-structures in $\EE$ (resp. strongly) trivializes the structure $S$.
	
An algebraic structure $S$ is self-chrystallographic when, for any category $\EE$, the category $S(\EE)$ is intensively crystallographic with respect to the structure $S$ itself.
\end{defi}

So, the algebraic structure $\Sigma$ trivializes the structure $S$ if and only if the canonical fully faithful embedding $Sub_1\EE \stackrel{\simeq}{\to} Sub_1\Sigma(\EE) \to S(\Sigma(\EE))$ is an equivalence of categories. It strongly trivializes it if and only if the canonical embedding $1 \into S(\Sigma(\EE))$ is an equivalence of categories. Since $S(\Sigma(\EE))=\Sigma(S(\EE))$, the \emph{trivialization} is clearly a symmetric relation between algebraic structures. On the other hand, an algebraic structure $S$ is self-chrystallographic if and only if the forgetful functor $U^S_{S(\EE)}: S(S(\EE)) \to S(E)$ is an equivalence of categories. The most basic examples are the following ones:
\\i) the structure of commutative monoid is self-chrystallographic;
\\ii) and the structure of abelian group as well;
\\iii) the structure of group and idempotent unitary magma strongly trivialize each other.
\proof
The two first points are straightforward. Let $(G,\circ,1)$ be a group and $*$ the binary law underlying the internal unitary magma structure in $\Gp$. By Eckmann-Hilton, we get $\circ=*$. From $x^2=x$, we get $x=1$, since this equality now holds inside a group.
\endproof
\begin{coro}\label{fund}
The fundamental group of a topological idempotent unitary magma is necessarily trivial.
More generally, when a structure $S$ is trivialized by the group structure, any topological $S$-structure produces a trivial fundamental group.
\end{coro}

\subsection{Subtractive categories}

The concept of subtractive category, introduced in \cite{ZJ}, is the categorical characterization of the pointed subtractive varieties in the sense of \cite{Urs}. A pointed category is subtractive when, given any reflexive relation $R$ on an object $X$, if $(0_X,1_X): X\into X\times X$ factorizes through $R$, then so does $(1_X,0_X): X\into X\times X$. It was shown in \cite{3} that on any object $X$ there is at most one structure of subtraction [$s(x,x)=0+s(x,0)=0$] and abelian group; so, the induced inclusion functors $\Ab\EE\into\Subt\EE\into \EE$ are fully faithfull, and we get:
\begin{coro}
Any subtractive category is crystallographic with respect to the structure of subtraction and abelian group.	
\end{coro}

\section{Further examples from the Mal'tsev setting}

\subsection{General Mal'sev setting}

A Mal'tsev structure is a set $X$ endowed with a ternary operation $p: X\times X\times X \to X$ satisfying the Mal'tsev identities $p(x,y,y)=x=p(y,y,x)$. We denote by $\Mal$ the category of  Mal'tsev structures. An affine structure is a Mal'tsev structure such that $p$ is associative ($p(x,y,p(z,u,v))=p(p(x,y,z),u,v)$) and commutative ($p(x,y,z)=p(z,y,x)$).  We denote by $\Aff$ the category of affine structures.

A category $\EE$ is a Mal'tsev one when any internal reflexive relation is actually an equivalence relation, \cite{CLP} and \cite{CPP}. Clearly any category $\Mal\EE$ of internal Mal'tsev  structures in $\EE$ is a Mal'tsev one.
Let us recall the following characterizations:
\begin{theo}\label{charmal}
The following conditions are equivalent:\\
1) $\EE$ is a Mal'tsev category;\\
2) any sub-reflexif graph of an internal groupoid in $\EE$ is an internal groupoid;\\
3) any fiber $Pt_Y{\EE}$ of the fibration of points is unital;\\
4) any fiber $Pt_Y{\EE}$ of the fibration of points is subtractive.
\end{theo}
The point 3) is grounded on the following definition: we denote by $Pt\mathbb{E}$ the category whose objects are the split epimorphisms in $\mathbb{E}$ with a given splitting and morphisms the commutative squares between these data; by $\P_{\mathbb E} :Pt\mathbb{E}\rightarrow \mathbb{E}$ we denote the functor associating its codomain with any split epimorphism. It is left exact and a fibration whose cartesian maps are the pullbacks of split epimorphisms; it is called the \emph{fibration of points}. The fiber above $Y$ is denoted by $Pt_Y\EE$. The points 2) and 3) come from  \cite{Bfib} and the point 4) from \cite{ZJ}.
The following observations are important for our purpose:
\begin{theo}\label{mal1}
Let $\EE$ be any Mal'tsev category. Then:\\
1) on any object $X$ there is at most one internal Mal'tsev structure which is necessarily an affine one;\\
2) on any internal reflexive graph in $\EE$, there is at most one internal category structure, which is necessarily a groupoid one;\\
3) any internal groupoid is necessarily an affine one.
\end{theo}
The two first points come from \cite{CPP}. A groupoid is affine when the associative Mal'tsev operation defined on the parallel arrows by $p(\phi,\chi,\psi)=\phi.\chi^{-1}.\psi$ is moreover commutative, and consequently produces an affine structure. The third point comes from \cite{Baff} where the \emph{affine groupoids} were introduced under the name of \emph{abelian groupoids}. From the point 1), we get immediately:
\begin{prop}
Given any Mal'tsev category $\EE$,\\
1) it is crystallographic with respect to the Mal'tsev structure and, a fortiori, with respect to the affine structure\\
2) any fiber $Pt_Y\EE$ is crystallographic with respect to the (abelian) group structure.	
\end{prop}

\subsection{Special Mal'sev settings}

Similarly to what happens in the unital setting and from the characterization Theorem  \ref{charmal}, there are two extremal cases of Mat'sev categories: when
i) any fiber $Pt_Y\EE$ is stiffly unital; and when
ii) any fiber $Pt_Y\EE$ is linear.

In the first case, the category $\EE$ is said to be a \emph{stiffly Mal'tsev} one, see \cite{Birk} from which we shall recall
 the three first points of the following characterization:
\begin{theo}\label{charstiff}
Let $\EE$ be any Mal'tsev category. The following conditions are equivalent:\\
1) the category $\EE$ is a stiffly Mal'tsev one;\\
2) the only abelian equivalence relations are the discrete ones;\\ 
3) any fiber $Pt_X\EE$ trivializes the group structure;\\
4) the only internal groupoids are the equivalence relations.
\end{theo}
The equivalence between 3) and 4) is shown below in Proposition \ref{truc}. Accordingly:
\begin{coro}
A stiffly Mal'tsev category $\EE$  trivializes the Mal'tsev structure. A stiffly Mal'tsev category $\EE$ strongly trivializes the group structure. 	
\end{coro}
\proof
In a Mal'tsev category, an object $X$ is affine if and only if the terminal map $\tau_X: X \to 1$ has an abelian kernel equivalence relation. So, our first assertion is a consequence of the point 2) of the previous characterization theorem which makes $\tau_X$ a monomorphism. An internal group object is always abelian in a Mal'tsev context; so, it is nothing but an affine object $X$ with a global element $e: 1 \into X$. The second assertion is then straightforward. 
\endproof

As for the second extremal case, let us recall  the following: 
\begin{theo}\label{mal2}
The following conditions are equivalent:\\
1) the category $\EE$ is a naturally Mal'tsev category;\\
2) on any reflexif graph in $\EE$ there is a unique structure of   internal groupoid;\\
3) any fiber $Pt_Y{\EE}$ is linear;\\
4) any fiber $Pt_Y{\EE}$ is additive.
\end{theo}
The equivalence between 1) and 2) is the ground result of \cite{john} where a naturally Mal'tsev category is defined as a category in which any object $X$ is endowed with a natural Mal'tsev structure. So, a naturally Mal'tsev category is a Mal'tsev one, and the natural Mal'tsev structure in question is then necessarily an affine one. Clearly the category $\Aff$ of affine structures and any category $\Aff\EE$ are naturally Mal'tsev ones. The equivalence between 1), 3) and 4) is in \cite{1}.
From that, we get immediately:
\begin{prop}
	1) Any naturally Mal'tsev category $\EE$ is intensively  crystallographic with respect to the Mal'tsev structure.\\
	 2) The affine structure is self-crystallographic. 	
\end{prop}

\section{Fiberwise extension of the notion of algebraic crystallography}\label{fw} 

From 2) in Theorem \ref{mal1} and in Theorem \ref{mal2}, we shall enlarge the set of examples of the crystallographic settings. For that, we shall think to the composition of arrows in a category or a groupoid as an algebraic structure on its underlying reflexive graph.

So, given any category $\EE$, let us denote by $RGh\EE$ the category of internal reflexif graphs in $\EE$:
$$\xymatrix@=8pt
{
	X_{\bullet}:\;\;\;\;\;  X_1 \ar@<2ex>[rrr]^>>>>>>>>{d_{1}^{X_{\bullet}}} \ar@<-2ex>[rrr]_>>>>>>>>{d_{0}^{X_{\bullet}}} &&&
	X_0   \ar[lll]|>>>>>>>>{s_0^{X_{\bullet}}}  
}
$$
which is finitely complete and by $U_0: RGh\EE \to \EE$ the left exact forgetful functor; it is a fibration whose cartesian map above $f_0: X_0\to Y_0$ is given by the inverse image along $f_0$. We denote by $RGh_X\EE$ the fibre above $X$; it has an initial object given by the discrete equivalence relation $\Delta_X$ on $X$ and a terminal one given by the indiscrete equivalence relation $\nabla_X$. Similarly we shall denote by $Cat\EE$ and $Grd\EE$ the categories of internal categories and groupoids whose objects are the 3-truncated simplicial objects in $\EE$ (including all the natural degeneracy maps which do not appear in the following diagram):
$$\xymatrix@=8pt
{
	X_{\bullet}:\;\;\;\;\; X_3\ar@<4ex>[rrr]^>>>>>>{d_{4}^{X_{\bullet}}} \ar@<2ex>[rrr]|>>>>>>{d_{3}^{X_{\bullet}}} \ar[rrr] \ar@<-2ex>[rrr]|>>>>>>{d_{1}^{X_{\bullet}}}\ar@<-4ex>[rrr]_>>>>>>{d_{0}^{X_{\bullet}}}	&&& X_2  \ar@<2ex>[rrr]^>>>>>>{d_{2}^{X_{\bullet}}} \ar[rrr]|>>>>>>{d_{1}^{X_{\bullet}}} \ar@<-2ex>[rrr]_>>>>>>{d_{0}^{X_{\bullet}}} &&& X_1 \ar@<2ex>[rrr]^{d_{1}^{X_{\bullet}}} \ar@<-2ex>[rrr]_{d_{0}^{X_{\bullet}}} &&&
	X_0   \ar[lll]|{s_0^{X_{\bullet}}}  
}
$$
where, in the first case, $X_2$ and $X_3$ are respectively obtained by the pullback of $d_0$ along $d_1$ and the pullback of $d_0$ along $d_2$ and where, in the second one, any commutative square is a pullback. The same inverse image process as above makes fibrations the following forgetful functors:
$$\xymatrix@=10pt
{
	Grd\EE \ar@{ >->}[rrr] \ar[dd]_{U_0} &&& Cat\EE \ar[dd]^{U_0}\ar[rrr] &&& RGh\EE\ar[dd]^{U_0} \\
	&&& &&& &&\\
	\EE \ar@{=}[rrr] &&& \EE \ar@{=}[rrr] &&& \EE 
}
$$
whose fibers, denoted respectively $Grd_X\EE$ and $Cat_X\EE$, have the same $\Delta_X$ and $\nabla_X$ as initial and terminal objects. Now the points 2) of Theorems \ref{mal1}, \ref{charstiff} and \ref{mal2} give rise to the following:
\begin{coro}
Let $\EE$ be a Mal'tsev category. Then any fiber $RGh_X\EE$:\\
1) is crystallographic with respect to the structure of affine groupoid;\\
2) trivializes 	the structure of groupoid if and only if $\EE$ is a stiffly Mal'tsev category;\\
3) is intensively crystallographic with respect to the structure of affine groupoid if and only if $\EE$ is a naturally Mal'tsev category.
\end{coro}
\proof
Only the point 2) demands some precisions: in a Mal'tsev category the subobjects of the terminal object $\nabla_X$ in the fiber $RGh_X\EE$ (namely the reflexive relations) are the equivalence relations.
\endproof

\section{Congruence hyperextensible and Gumm categories}

\subsection{Congruence modular varieties and Gumm categories}

A congruence modular variety is a variety in which the modular formula holds for congruences:
$$(T\vee S) \wedge R=T\vee (S\wedge R), \rm{\; for \; any \; triple\;} (T,S,R)\; \rm{\; such \; that:\;} T\subset R$$

In \cite{Gumm}, they were characterized in geometric terms by the validity of the \emph{Shifting Lemma}:
given any triple of equivalence relations $(T,S,R)$ such that $R\cap S\subset T$, the following left hand side situation implies the right hand side one:
$$ \xymatrix@=15pt{
	x \ar[r]^{S} \ar@(l,l)[d]_T \ar[d]_R & y \ar[d]^R &&  y \ar@(r,r)[d]^T\\
	x' \ar[r]_{S} & y' &&   y'
}
$$
One of the main interest of the Shifting lemma is that it is freed of any condition involving finite colimits. So, thanks to the Yoneda embedding, it  keeps a meaning in any finitely complete category $\EE$. This led to the notion of Gumm category \cite{BGGGum}. Any regular Mal'tsev category is a Gumm one. Let us recall, from the same article, the following:
\begin{prop}
Let $\EE$ be any Gumm category. Then:\\
1) on any object $X$ there is at most one affine structure;\\
2) on any reflexive graph, there is at most one structure of category.	
\end{prop}
The first point was already noticed by Gumm in the congruence modular varieties. On the contrary to what happens in Mal'tsev ones, in Gumm categories there are internal categories which are not groupoids. Now we get:
\begin{coro}
Let $\EE$ be any Gumm category:\\
1) it is chrytallographic with respect to the affine structure;\\
2)  any fiber $RGh_X\EE$ is chrytallographic with respect to the structure of category.		
\end{coro}

\subsection{Conruence hyperextensible categories}

On the model of Mal'tsev ones, Gumm categories have a characterization through the fibration of points \cite{BGum}:
\begin{theo}
	Given a category $\EE$, the following conditions are equivalent:\\
	1) $\EE$ is Gumm category;\\
	2) any fiber $Pt_Y\EE$ is congruence hyperextensible.
\end{theo}

\begin{defi}
	A pointed category $\EE$ is said to be congruence hyperextensible when given any punctual span in $\EE$ (= any commutative diagram of split epimorphisms):
	$$\xymatrix@C=2,5pc@R=1,5pc{ W \ar[d]_{f}\ar@<1ex>[r]^g & Y\ar[d]_{\tau_Y}\ar@{ >->}[l]^{t}\\
		X\ar@<-1ex>@{ >->}[u]_{s} \ar@<-1ex>[r]_{\tau_X} & 1\ar@<-1ex>[u]_{0_Y} \ar[l]_{0_X} }$$
and any equivalence relation $T$ on $W$ with $R[f]\cap R[g]\subset T$, we get:\\ $R[f]\cap g^{-1}(t^{-1}(T))\subset T$.  
\end{defi}
This means that the following situation holds in $W$:
$$ \xymatrix@=15pt{
	tg(x) \ar[r]^{R[g]} \ar@(l,l)[d]_T \ar[d]_{R[f]} & x \ar[d]^{R[f]} &&  x \ar@(r,r)[d]^T\\
	tg(x') \ar[r]_{R[g]} & x' &&   x'
}
$$
Any J\'onsson-Tarski variety $\VV$ or regular unital category $\EE$ is congruence hyperextensible. It is a fortiori the case for any regular pointed protomodular category \cite{BB}, and thus for any Slomi\'nsky variety \cite{slo}. Here are some examples of varieties which are outside of these classes:
\begin{prop}\label{1}
	Let $\VV$ be a pointed variety with $2k+1$ ternary terms $p_i$  such that:
	\begin{align*}
	p_1(a,0,0)=a \;\; ,&  \;\;\;  p_i(a,0,a)=a  \; \;\;  2\leq i \leq 2k \;\;, & p_{2k+1}(0,0,a)=a \\
	p_{2i-1}(a,a,b)=p_{2i}(a,a,b) &  \;\;\;\;  1\leq i\leq k & \\
	p_{2i}(a,b,b)=p_{2i+1}(a,b,b) &  \;\;\;\;   1\leq i\leq k & \\ 
	& \;\;\;\;\;\;\;  p_i(b,b,b)=b,  \;  & \;  1\leq i \leq 2k+1  & 
	\end{align*}
	Then the variety $\VV$ is congruence hyperextensible but not a J\'onsson-Tarski one.
\end{prop}
\proof
Starting with a punctual relation in $\VV$ as above
and with an equivalence relation $T$ on $W$ such that $(0Wb)T(0Wc)$.
First let us show that:\\ $p_2(aWb,0Wc,aWc)Tp_{2k}(aWb,0Wc,aWc)$. For all $i$ with  $ 1\leq i \leq k$, we get:\\
$p_{2i}(aWb,0Wc,aWc)=p_{2i}(a,0,a)Wp_{2i}(b,c,c)=p_{2i+1}(a,0,a)Wp_{2i+1}(b,c,c)=$\\
$p_{2i+1}(aWb,0Wc,aWc)$. So: $p_{2i}(aWb,0Wc,aWc)Tp_{2i+1}(aWb,0Wb,aWc)$.\\ 
Then: $p_{2i+1}(aWb,0Wb,aWc)=p_{2i+1}(a,0,a)Wp_{2i+1}(b,b,c)$\\
$=p_{2i+2}(a,0,a)Wp_{2i+2}(b,b,c)=p_{2i+2}(aWb,0Wb,aWc)$, and finally:\\$p_{2i+2}(aWb,0Wb,aWc)Tp_{2i+2}(aWb,0Wc,aWc)$.

We get also: $aWb=p_1(a,0,0)Wp_1(b,b,b)=p_1(aWb,0Wb,0Wb)$.
Whence:\\ $(aWb)Tp_1(aWb,0Wb,0Wc)$ $=(aWb)T(aWp_1(b,b,c))=(aWb)T(aWp_2(b,c,c))$\\
$=(aWb)T(p_2(a,0,a)Wp_2(b,c,c))=(aWb)Tp_2(aWb,0Wc,aWc)$.

Finally: $p_{2k}(aWb,0Wc,aWc)=p_{2k}(a,0,a)Wp_{2k}(b,c,c)$\\$=p_{2k+1}(0,0,a)Wp_{2k+1}(b,c,c)=p_{2k+1}(0Wb,0Wc,aWc)$. Whence finally:\\
$p_{2k}(aWb,0Wc,aWc)Tp_{2k+1}(0Wc,0Wc,aWc)=p_{2k}(aWb,0Wc,aWc)T(aWc)$.

The last assertion is proved when $n=3$ in \cite{2}.
\endproof
We shall say that a congruence hyperextensible variety of the previous kind is a CHyper variety of type $2k+1$. \emph{Any CHyper variety of type $2k+1$ is of type $2k+3$}:
\proof
Let $\VV$ be a CHyper variety of type $2k+1$. Let us set: $$\bar p_i=p_i, \; \forall i \; 1\leq i\leq 2j; \; \bar p_{2j+1}=p_{2j}=\bar p_{2j+2}; \; \bar p_i=p_{i-2}, \;\; \forall i \;\; 2j+3\leq i\leq 2k+1$$
These $\bar p_i$ make the variety  $\VV$ a  CHyper variety  of type $2k+3$.
\endproof
Recall now the main observation from \cite{2}:
\begin{prop}\label{abob}
	Let $\EE$ be a congruence hyperextensible category. Any subtraction $s$ on an object $X$ is the difference mapping associated with an internal group structure which is necessarily abelian. On any object $X$, there is at most one subtraction $s$; any morphism $f:X\to Y$ between objects with subtraction is a subtraction homomorphism. 
\end{prop}
Whence the following result which, as explained in the introduction, led to the notion: 
\begin{coro}
Any congruence hyperextensible category is crystallographic with respect to structure of subtraction and (abelian) group. Given a Gumm category $\EE$, any fiber $Pt_Y\EE$ is crystallographic with respect to the (abelian) group structure.	
\end{coro}

\subsection{Congruence distributive varieties and categories}

It is well known that any congruence distributive variety in which the following congruence formula holds:
$$(T \wedge R)\vee (T \wedge S)=T \wedge (R \vee S)$$
is congruence modular. Let us recall that a Mal'tsev variety $\VV$ is congruence distributive if and only if its Mal'tsev term $p$ satisfies the  Pixley axiom $p(x,y,x)=x$ \cite{Pix}. This is the case in particular for the category $\Heyt$ of Heyting algebras. In this section, we shall show that congruence distributive varieties and categories trivialize the group structure.
\begin{defi}
	A category $\EE$ is said to be weakly congruence distributive when  the following implication holds:
	$$(T \wedge R=\Delta_X \; {\; \rm and \;} \; T \wedge S=\Delta_X) \Rightarrow T \wedge (R \vee S)=\Delta_X$$
	whenever the supremum $R \vee S$ of the pair $(R,S)$
	of equivalence relations does exist.
\end{defi}
\begin{prop}
Given any weakly congruence distributive category $\EE$, any fiber $Pt_X\EE$ trivializes the group structure; the category $\EE$ trivializes the associative Mal'tsev structure.
\end{prop}
\proof
In any category $\EE$, given any pair $(X,Y)$ of objects, we get:\\ $R[p_X]\vee R[p_Y]=\nabla_{X\times Y}$, since we get it in $Set$.\\
1) Let us begin by showing that $\EE$ trivialize the group structure. So, let us denote by $m: X\times X \to X$ the binary operation associated with an internal monoid structure. It is a group structure if and only if the following square is a pullback:
$$\xymatrix@C=1,5pc@R=1,5pc{ X\times X \ar[d]_{p_0}\ar[rr]^m && X\ar[d]\\
	X\ar@<-1ex>@{ >.>}[u]_{(1_X,()^{-1})} \ar[rr]_{\tau_X}
	 && 1\ar@{>.>}@<-1ex>[u]_{e}  }$$
 Then the unit $e$ produces the left hand side section.
 Now, observe that both $(m,p_0)$ and  $(m,p_1)$ are monomorphisms; this means that $$R[m] \wedge R[p_0]=\Delta_{X\times X} \; {\; \rm and \;} \; R[m] \wedge R[p_1]=\Delta_{X\times X}$$ Whence: $\Delta_{X\times X}=R[m] \wedge (R[p_0] \vee R[p_1])=R[m]\vee \nabla_{X\times X}=R[m]$. Accordingly the morphism $m$ is a monomorphism; being split (by $e\times 1_X$), it is an isomorphism. Coming back to the pullback above, the upward square is a pullback, and $\tau_X$ is an isomorphism as well. Accordingly, $X$ is isomorphic to the terminal object.\\
 2) We can then reproduce this proof in any slice category $\EE/X$ or equivalently in any fiber $Pt_X\EE$, since, similarly to 1), given any pullback in $\EE$:
 $$\xymatrix@C=2,5pc@R=1,5pc{ X\times_ZY \ar[d]_{p_X} \ar@{.>}[dr]_{\gamma} \ar[r]^{p_Y} & Y\ar[d]^{g}\\
 	X \ar[r]_{f} & Z }$$
 we get $R[p_X]\vee R[p_Y]=R[\gamma]$ in $\EE$: if $S$ is any equivalence relation on $X\times_ZY$ in $\EE$ such that $R[p_X]\subset S$ and $R[p_Y]\subset S$, then the same inclusions hold for $S\wedge R[\gamma]$ which now lies in $\EE/Y$. So, we get $R[\gamma]=S \wedge R[\gamma]$ in the category $\EE/Y$, whence $R[\gamma]\subset S$ in $\EE$.\\
 3) Now let $p: X\times X \times X \to X$ be any internal associative Mal'tsev structure. We then get an internal group structure  on the object $(p_0^X,s_0^X): X\times X \splito X$ of the fiber $Pt_X\EE$, by setting $(x,y)*(x,z)=(x,p(y,x,z))$. So, $p_0^X$ is an isomorphism; accordingly $\tau_X: X\to 1$ is a monomorphism and $X$ a subobject of $1$.
\endproof 
It is a fortiori the case for any congruence distributive variety or category. In the Mal'tsev context, we have the converse with the following precisions:
\begin{prop}
Given any regular Mal'tsev category $\EE$, the following conditions are equivalent:\\
1) $\EE$ is weakly congruence distributive;\\
2) any fiber $Pt_X\EE$ trivializes the group structure;\\ 3) the category $\EE$ is a stiffly Mal'tsev one.\\
Given any exact Mal'tsev category $\EE$, the following conditions are equivalent:\\
1) $\EE$ is  congruence distributive;\\
2) any fiber $Pt_X\EE$ trivializes the group structure;\\ 3) the category $\EE$ is a stiffly Mal'tsev one.	
\end{prop}
\proof
The first equivalence comes from Lemma 2.9.10 and Proposition 1.11.34 in \cite{BB}, while the second one is Theorem 2.9.11.
\endproof

Accordingly, a Mal'tsev variety is a stiffly Mal'tsev category if and only if it satisfies the Pixley axiom.  Finally let us observe that more generally we get the following:
\begin{prop}\label{truc}
Given any category $\EE$, the following conditions are equivalent:\\
1) 	any fiber $Pt_X\EE$ trivializes the group structure;\\
2) any internal groupoid is an equivalence relation.\\
Such a category $\EE$ trivializes the associative Mal'tsev structure.
\end{prop}
\proof
 We shall work in the fiber $Grd_X\EE$ which is a protomodular category \cite{BB}; accordingly, in this fiber, pullbacks reflect monomorphisms. Let $X_{\bullet}$ be any groupoid in this fiber and consider the following pullback in $Grd_X\EE$:
 $$\xymatrix@C=1,5pc@R=1,5pc{ EndX_{\bullet} \ar[d]\ar@{ >->}[rr] && X_{\bullet} \ar[d]\\
 	\Delta_X\ar@<1ex>@{ >->}[u]_{} \ar@{ >->}[rr]
 	&& \nabla_X  }$$
 where $EndX_{\bullet}$ is the sub-groupoid of the endomorphisms of $X_{\bullet}$. It coincides with a group structure in the fiber $Pt_X\EE$. So, when 1) holds, the left hand side vertical map is an isomorphism, and consequently the right hand side one is a monomorphism, and thus $X_{\bullet}$ is an equivalence relation. Conversely any  group structure in the fiber $Pt_X\EE$ produces a groupoid $X_{\bullet}$ which coincides with its endsome $EndX_{\bullet}$; in other words the upper horizontal arrow in the previous diagram is an isomorphism. So, if we suppose 2), namely that the right hand side vertical map is a monomorphism, so is the left hand side one. Being split, this same map is an isomorphism, and the group in question is trivial in $Pt_X\EE$.
 The last assertion is checked as in the previous proposition.
\endproof
In this context, we then get another possible level of trivialization:
\begin{defi}
A category $\EE$ weakly trivializes a structure $S$ when the only $S$-structures in $\EE$ are given by \emph{some class} of subobjects of $1$.	
\end{defi}
\begin{prop}
When $\EE$ is a category such that 	any fiber $Pt_X\EE$ trivializes the group structure, then any fiber $RGh_X\EE$ weakly trivializes the structure of groupoid.
\end{prop}

\section{First principles and questionings about the algebraic crystallography}

\subsection{Commutativity}

The three first examples of crystallographic context (unital, additive, Mal'tsev) are actually dealing with commutative structures. The reason is simple. Any of the structures $S$ in question (monoid, group, Mal'tsev structure) is endowed with a duality operator making the following diagram commute:
$$\xymatrix@=6pt
{
	S(\Set) \ar[rrrr]^{(\;)^{op}} \ar[ddrr]_{U^S} &&&&  S(\Set) \ar[ddll]^{U^S}\\
	&&& &&& \\
 && \Set 
}
$$
The uniqueness involved in the definition of a crystallographic context implies that on an object $X$ any structure coincides with its dual.

\subsection{Questions}

Of course, a first natural question is: is it possible to characterize in some way the structures which admit a crystallographic context or an intensive crystallographic context? Concerning the classical (=non-fiberwise) algebraic structures, we can ask when the first situation implies the second one.

Another way to think to this question would be: are the commutative monoid, abelian group and affine structure the only ones to be  self-crystallogra- phic structures?

\subsection{Trivializing chrystallographic context}

Certainly, as we already noticed, trivializing contexts are not so scarce. Let us add the following very special one, recalling from \cite{abbott} the following:

\begin{defi}
An implication algebra is a triple $(X,\trio,1)$ of a set, a binary operation and a constant such that:
\begin{align*}
x\trio x=1 \;\; (1) & & (x\trio y)\trio x=x \;\; (2)\\
(x\trio y)\trio y=(y\trio x)\trio x \;\; (3) & & x\trio (y \trio z)=y \trio (x\trio z) \;\; (4)
\end{align*}	
\end{defi}
\begin{lemma}
In any impication algebra we get:
$$1\trio x=x \;\;\; {\rm and} \;\;\; x\trio 1=1$$
Accordingly any impication  algebra has an underlying opsubtraction.
\end{lemma}
\proof
First set $y=x$ in (2). Then set $x=y$ in the same (2).
\endproof
Let us denote by $\Imp$ the category of implication algebras.
\begin{prop}
	The category $\Imp$ trivializes the structure of implication algebra. In other words, the structure of implication algebra is self-trivializing.
\end{prop}
\proof
Let $*$ give an internal implication algebra structure on $(X,\trio,1)$ in $\Imp$. It must be a $\Imp$-homomorphism, that is to satisfy:
$$(x\trio x')*(y\trio y')=(x*y)\trio(x'*y')$$
Setting $x=x'=y$, we get:
 $y\trio y'=y*y'$. So we must have:
$$(x\trio x')\trio(y\trio y')=(x\trio y)\trio(x'\trio y')$$
Setting $x=x'=y'$ and $y=1$ , we get: $y'=1$; so, $(X,\trio,1)$ is the terminal object in $\Imp$.
\endproof 
According to the proofs of the previous lemma and proposition, the same result holds for the structure of \emph{implicative opsubtraction} defined by $s(x,x)=0, s(0,x)=x$ and $s(x,0)=0$ or, by duality,  for the structure of \emph{opimplicative subtraction}. Let $(X,p)$ be any Mal'tsev structure such that $p$ satisfies the Pixley axiom; then, given any element $x_0\in X$, the binary operation defined by $s_{x_0}(y,z)=p(y,z,x_0)$ produces an opimplicative subtraction on $X$ with $x_0$ as constant. Incidentally, we get the following:
\begin{prop}
The structure of opimplicative subtraction (resp. implicative opsubtraction) trivializes the structure of unitary magma. The fundamental group of a topological 	opimplicative subtraction (resp. implicative opsubtraction) is trivial. Accordingly, the fundamental group, at any point, of a topological Mal'tsev algebra satisfying the Pixley axiom is trivial.
\end{prop}
\proof
Let $(X,*)$ be a unitary magma, and $s$ the internal opimplicative subtraction on it. We have: $s(x*x',y*y')=s(x,y)*s(x',y')$. Taking $x'=y=1$, we get: $s(x,y')=x*1=x$. Then take $x=y'$; so, $x=1$. The second assertion comes from Corollary \ref{fund}. A topological Mal'tsev algebra is any object of a category $\mathbb{T}$-$Top$ of internal $\mathbb{T}$-algebras in the category $Top$ of topological spaces, where $\mathbb{T}$ is a Mal'tsev variety, see \cite{JP} for instance. The last assertion is then straightforward from the opimplicative subtraction $s_{x_0}$. 
\endproof

\subsection{Paradoxical aspect of the notion}

On the one hand, this kind of relationship \emph{"context vs structure"} has a classical positive aspect on the side of the context:\\
\emph{in this given context, any object $X$ has at most one structure}.

On the other hand, on the side of the structure, it could be thought as a kind of photographic negative:\\\emph{on any object in the context, there is no more than one structure of this type}.

So, emerge a paradoxical question: could it be possible to extract some  interesting (positive) information about a structure from the contexts in which this structure becomes so punctually scarce? Or, for example and more precisely, what are we exactly learning about the notion of (abelian) group by knowing that congruence hyperextensive categories produce a crystallographic context with respect to it?

\section{Some very large abelian and naturally Mal'tsev categories}

As said in the introduction, in this last section, we shall produce some constructions which seem to offset the punctual scarcity of a structure in a chrystallographic context by a kind of multiplication of this structure inside this context; first in the pointed case of abelian groups, then in the non-pointed case of affine structures.

\subsection{A very large abelian category}\label{alab}

Consider the congruence hyperextensible variety $Hex_3$ defined by the only following three ternary terms and the following equations:
	\begin{align*}
p_1(a,0,0)=a, & &  p_2(a,0,a)=a, & & p_{3}(0,0,a)=a, \\
p_{1}(a,a,b)=p_{2}(a,a,b), & & p_{2}(a,b,b)=p_{3}(a,b,b), & & p_i(b,b,b)=b, \;\;  1\leq i \leq 3. 
\end{align*}
We get a fully faithful functor $h: \Gp\to Hex_3$ defined in the following way: starting from a group $(G,.)$, construct a $Hex_3$-algebra on the set $G$ by setting: $$p_1(x,y,z)=x.y^{-1}.z, \;\; \ {\rm and} \;\; p_2(x,y,z)=z=p_3(x,y,z)$$ 
By restriction we get a fully faithful embedding $h:\Ab \to \Ab(Hex_3)$.

Now consider any field $K$ whose characteristic is not $2$. We get a faithful functor $w_{K}: K$-$Vect \to \Ab(Hex_3)$: starting from a $K$-vector space $V$, construct a $Hex_3$-algebra on the set $V$ by setting: 
$$\bar p_1(x,y,z)=x+ \frac{-y+z}{2}, \;\; \bar p_2(x,y,z)=\frac{x+z}{2}, \;\; {\rm and} \;\;\bar p_3(x,y,z)=\frac{x-y}{2}+z$$ 
This algebra is made abelian by the $Hex_3$-homomorphism $+: V\times V \to V$.

So, the category $\Ab Hex_3$ of abelian objects in $Hex_3$ appears:\\ 
i) to fully faithfully embed the category $\Ab$ of abelian groups by the functor $h: \Ab \rightarrow \Ab Hex_3$\\ and in an independant way\\
ii) to faithfully contain any category $K$-$\Vect$ of $K$-vector spaces, provided that the field $K$ is not of characteristic $2$, by the functor $w_{K}: K$-$\Vect\rightarrow \Ab Hex_3$.\\
So, the underlying abelian group  $(V,+)$ of a $K$-vector space is represented by two distinct  objects, $h(V,+)$ and $w_K(V,+)$, in the abelian category $\Ab Hex_3$.

\subsection{A very large exact naturally Mal'tsev category}

We are now going to extend the previous kind of situation to a non-pointed context.

\medskip 

Consider the congruence modular variety $CM_3$ defined by the only following three ternary terms and the following equations:
\begin{align*}
p_1(a,b,b)=a, & &  p_2(a,b,a)=a, & & p_{3}(b,b,a)=a, \\
p_{1}(a,a,b)=p_{2}(a,a,b), & & & & p_{2}(a,b,b)=p_{3}(a,b,b). 
\end{align*}
We get a fully faithful functor $m: \Mal \to CM_3$ defined in the following way: starting from an algebra $(X,p)$ in $Mal$, construct a $CM_3$-algebra on the set $X$ by setting:
$$p_1=p \;\;\;\;\;\; {\rm and} \;\;\;\;\;\; p_2(x,y,z)=z=p_3(x,y,z)$$
By restriction, we get a fully faithful functor $m:\Aff \to \Aff(CM_3)$.

Now consider any field $K$ whose characteristic is not $2$. We get a faithful functor $a_{K}: K$-$\Aff \to \Aff(CM_3)$: starting from a $K$-affine space $X$, construct an affine $CM_3$-algebra on the set $X$ by setting: $$\bar p_1(x,y,z)=\beta(\dot x+ \frac{-\dot y+ \dot z}{2}), \;\; \bar p_2(x,y,z)=\beta(\frac{\dot x+ \dot z}{2}), \;\; \bar p_3(x,y,z)=\beta(\frac{\dot x-\dot y}{2}+z)$$
where $\beta$ is the barycentric mapping which is actually a $K$-$\Aff$-homomorphism. The affine structure $X\times  X \times X \to X$ on the algebra $a_K(X)$ in the variety $CM_3$ is then produced by the $CM_3$-homomorphism $p(x,y,z)=\beta(\dot x -\dot y+\dot z)$. 

So, the category $\Aff CM_3$ of affine structures in the variety $CM_3$ appears:\\ 
i) to fully faithfully embed the category $\Aff$ of affine structures by the functor $m: \Aff \rightarrow \Aff CM_3$\\ and in an independant way\\
ii) to faithfully contain any category $K$-$\Aff$ of $K$-affine spaces, provided that the field $K$ is not of characteristic $2$, by the functor $a_{K}: K$-$\Aff\rightarrow \Aff CM_3$.\\
So, the affine structure $(X,\beta)$ underlying a $K$-affine space $(X,\beta)$ is represented by two distinct objects, $m(X,\beta)$ and $a_K(X,\beta)$, in the naturally Mal'tsev category $\Aff CM_3$.

\noindent Keywords: Internal structures, Mal'tsev, congruence modular and congruence distributive varieties, unital, Mal'tsev and Gumm categories, context vs structure.\\
\noindent Amsclass: 08A05,08B05,08B10,18C10,18C40,18E13.

\vspace{3mm}
\noindent Univ. Littoral C\^ote d'Opale, UR 2597, LMPA,\\
Laboratoire de Math\'ematiques Pures et Appliqu\'ees Joseph Liouville,\\
F-62100 Calais, France. bourn@univ-littoral.fr
\end{document}